  %

 \input amstex.tex
 \input amsppt.sty
\magnification=\magstep1
\def\nava{\leavevmode\hbox{,\kern-0.09em,}}
\def\latinicacm{\def\DJ{\leavevmode\setbox0=\hbox{D}\kern0pt\rlap
 {\kern.04em\raise.188\ht0\hbox{-}}D}\def\dj{\leavevmode
 \setbox0=\hbox{d}\kern0pt\rlap{\kern.215em\raise.46\ht0\hbox{-}}d}\def\SH{\v
 S}\def\sh{\v s}\def\CJ{\' C}\def\cj{\' c}\def\CH{\v C}\def\ch{\v c}\def\Z{\v
 Z}\def\z{\v z}\def\DZ{D{\z}}\def\dz{d{\z}}}
\latinicacm

\def\bsk{\bigskip}
\def\cent{\centerline}
\def\msk{\medskip}
\def\ssk{\smallskip}

\def\ltextindent#1{\hbox to \hangindent{#1\hss}\ignorespaces }%
\newbox\uzezabox
\def\uzeza#1{\setbox\uzezabox=\hbox{#1}\advance\parindent by\wd\uzezabox}
\font\caps=cmcsc10

\def\mcirc{\mathbin{\raise0.8pt\hbox{$\ssize\circ$}}}

\def\iADS#1#2#3{\abovedisplayskip#1pt plus#2pt minus#3pt}
\def\iBDS#1#2#3{\belowdisplayskip#1pt plus#2pt minus#3pt}
\def\iADSS#1#2#3{\abovedisplayshortskip#1pt plus#2pt minus#3pt}
\def\iBDSS#1#2#3{\belowdisplayshortskip#1pt plus#2pt minus#3pt}
\catcode`\@=11
\def\foliofont@{\tenrm}
\def\logo@{}
\font@\tencyr=wncyr10
\font@\tencyrb=wncyb10
\font@\tencyri=wncyi10
\font@\eightcyr=wncyr8
\font@\eightcyrb=wncyb8
\font@\eightcyri=wncyi8
\font@\ninecyr=wncyr9
\font@\ninecyrb=wncyb9
\font@\ninecyri=wncyi9
\catcode`\@=13
\def\cirilicawn{\def\1{{\kern0pt}}\def\DJ{\char'006}\def\dj{\char'016}\def
 \SH{\char'130}\def\sh{\char'170}\def\CJ{\char'007}\def\cj{\char'017}\def
 \CH{\char'121}\def\ch{\char'161}\def\Z{\char'021}\def\DZ{\char'002}\def
 \dz{\char'012}\def\z{\char'031}
 \def\JO{\char'023}\def\jo{\char'033}\def\JU{\char'020}\def\ju{\char'030}\def
 \JI{\char'022}\def\ji{\char'032}\def\JA{\char'027}\def\ja{\char'037}\def
 \TV{\char'137}\def\tv{\char'177}\def\M{\char'136}\def\m{\char'176}\def
 \EE{\char'003}\def\ee{\char'013}\def\NO{\char'175}\def\rnava{\char'074}\def
 \rnavb{\char'076}}
\addto\tenpoint{}
\addto\eightpoint{}
\addto\tenpoint{}
\addto\eightpoint{}
\def\DisplaysX{\iADS{12}{3}{9}\iADSS{0}{3}{0}%
\iBDS{12}{3}{9}\iBDSS{7}{3}{4}}%
\def\DisplaysVIII{\iADS{9.6}{2.4}{7.2}\iADSS{0}{2.4}{0}%
\iBDS{9.6}{2.4}{7.2}\iBDSS{5.6}{2.4}{3.2}}%
\def\DisplaysVII{\iADS{8.4}{2.1}{6.3}\iADSS{0}{2.1}{0}%
\iBDS{8.4}{2.1}{6.3}\iBDSS{4.9}{2.1}{2.8}}%
\def\DisplaysV{\iADS{6}{1.5}{4.5}\iADSS{0}{1.5}{0}%
\iBDS{6}{1.5}{4.5}\iBDSS{3.5}{1.5}{2}}%
\def\strat{\hbox{\vrule height10.5pt depth0pt width0pt}}
\def\ESE{\eightpoint\strat\eightrm }
\tenpoint\DisplaysX
\addto\eightpoint{\normalbaselineskip9pt\normalbaselines\DisplaysVIII}
\hsize=6.5truein
\vsize=8.0truein
\hoffset=0 truept
\voffset=0 truept
\lefthyphenmin=2
\righthyphenmin=2
\headline={\hfil}
\footline={{\hfil\tenrm\folio\hfil}}
\output={\plainoutput}



\font\ninerm =cmr9
\font\ninei  =cmmi9
\font\ninesy =cmsy9
\font\ninebf =cmbx9
\font\ninett =cmtt9
\font\nineit =cmti9
\font\ninesl =cmsl9
\font\sevensl=cmsl8 scaled 913
\font\seventt=cmtt8 scaled 913

\newskip\ttglue

\def\ninepoint{\def\rm{\fam0\ninerm} 
  \textfont0=\ninerm  \scriptfont0=\sixrm      \scriptscriptfont0=\fiverm
  \textfont1=\ninei   \scriptfont1=\sixi       \scriptscriptfont1=\fivei
  \textfont2=\ninesy  \scriptfont2=\sixsy      \scriptscriptfont2=\fivesy
  \textfont3=\tenex   \scriptfont3=\tenex      \scriptscriptfont3=\tenex
  \textfont\itfam=\nineit  \def\it{\fam\itfam\nineit}%
  \textfont\slfam=\ninesl  \def\sl{\fam\slfam\ninesl}%
  \textfont\ttfam=\ninett  \def\tt{\fam\ttfam\ninett}%
  \textfont\bffam=\ninebf  \def\bf{\fam\bffam\ninebf}%
  \scriptfont\bffam=\sevenbf \scriptscriptfont\bffam=\fivebf
  \tt \ttglue=.5em plus .25em minus .15em
  \normalbaselineskip=11pt
  \setbox\strutbox=\hbox{\vrule height8pt depth3pt width0pt}%
  \let\sc=\sevenrm     \let\big=\ninebig    \normalbaselines\rm}

 \def\tenbig#1{{\hbox{$\left#1\vbox to8.5pt{}\right.\n@space$}}}
 \def\ninebig#1{{\hbox{$\textfont0=\tenrm\textfont2=\tensy
                       \left#1\vbox to7.25pt{}\right.\n@space$}}}
 \def\eightbig#1{{\hbox{$\textfont0=\ninerm\textfont2=\ninesy
                       \left#1\vbox to6.5pt{}\right.\n@space$}}}

\def\sevenpoint{\def\rm{\fam0\sevenrm} 
  \textfont0=\sevenrm  \scriptfont0=\sixrm      \scriptscriptfont0=\fiverm
  \textfont1=\seveni   \scriptfont1=\sixi       \scriptscriptfont1=\fivei
  \textfont2=\sevensy  \scriptfont2=\sixsy      \scriptscriptfont2=\fivesy
  \textfont3=\tenex    \scriptfont3=\tenex      \scriptscriptfont3=\tenex
  \textfont\itfam=\sevenit  \def\it{\fam\itfam\sevenit}%
  \textfont\slfam=\sevensl  \def\sl{\fam\slfam\sevensl}%
  \textfont\ttfam=\seventt  \def\tt{\fam\ttfam\seventt}%
  \textfont\bffam=\sevenbf  \def\bf{\fam\bffam\sevenbf}%
  \scriptfont\bffam=\sevenbf \scriptscriptfont\bffam=\fivebf
  \tt \ttglue=.4em plus .15em minus .10em
  \normalbaselineskip=8pt
  \setbox\strutbox=\hbox{\vrule height8pt depth3pt width0pt}%
  \let\sc=\sixnrm     \let\big=\sevenbig    \normalbaselines\rm}
 \topmatter
  \tolerance=500
  \emergencystretch=0pt
  \nopagenumbers
\hsize=12.5cm
\vsize=17cm
\hoffset=1cm
\voffset=2cm
  \parindent=25pt
  \parskip\smallskipamount
  \MultlineGap{1em}
  \pageno=1
  \def\rightheadline{{\hfil{\eightpoint\eightrm\rm
  On Kurepa's problems in number theory
  }\hfil\tenrm\folio}}
  \def\leftheadline{{\tenrm\folio\hfil{\eightpoint\eightrm \rm
  A. Ivi\'c, \v{Z}. Mijajlovi\'c
  }\hfil}}
  \def\emptyheadline{\hfil}
  \headline{\ifnum\pageno=1 \emptyheadline\else
  \ifodd\pageno \rightheadline \else \leftheadline\fi\fi}
  \TagsOnRight
  \DisplaysVII
  \addto\eightpoint{\DisplaysV}
  \def\tompublikasiona{57 (71)}%
  \def\godistepublikasiona{1995}%
  \def\stranicepublikasiona{19--28}%
  \setbox1=\hbox{\hss{\eightrm  {\DJ}uro Kurepa memorial volume}\hss}  
  \setbox0=\hbox{\sevenrm PUBLICATIONS DE L'INSTITUT MATH\'EMATIQUE}
  \leftline{\vbox{\baselineskip=8pt
  \copy0\hbox to \wd0
  {\hss\sevenrm  Nouvelle s\'erie, tome {\tompublikasiona},
  {\godistepublikasiona}, {\stranicepublikasiona}\hss}
  \hskip 5.5mm\box1                        
  }}
  \vskip105pt\baselineskip=12pt
  \vskip -4mm                              

  \def\GF{\text{\rm GF}}
  \def\rest{\text{\rm rest}}
  \def\mmod{\text{\rm mod}\,}
  \def\Re{\text{\rm Re}}

  \cent{\bf        ON KUREPA'S PROBLEMS IN NUMBER THEORY
  }
  \smallskip
  \cent{\it      Dedicated to the memory of Prof\. {\DJ}uro Kurepa
  }
  \bsk
   \cent{\bf      A. Ivi\'c, \v{Z}. Mijajlovi\'c
    }
  \bsk\msk

  \begingroup
  \eightpoint

  {\bf Abstract}.
       We discuss some problems in number theory posed by \DJ uro Kurepa,
  including the so-called left factorial hypothesis that an odd prime
  $p$ does not divide $0! + 1! + \cdots + (p-1)!$.
  \par
  \endgroup
  \bsk


   \footnote""{\ESE{\it AMS Subject Classification\/} (1991):
   Primary 11\,A\,05}

  \centerline{\bf 1. Introduction}
  \smallskip

  \DJ. Kurepa posed several problems in number theory that drew attention of
  many workers in number theory.
  Certainly, the most known of his problems is the so called left factorial
  hypothesis, which is still an open problem. However, Kurepa asked
  several other questions that are less known, but we think that they are
  interesting as well. The aim of this paper is to review some of these
  problems, and to present some of the known results  concerning them.

      We shall assume the following notation. We shall denote by $N$
  the set of natural numbers (nonnegative integers), $N^+$ denotes positive
  integers, while ${\bold Z}_n$ denotes the ring of integers modulo $n$.
  The greatest common divisor of integers $a$ and $b$ is denoted by $(a,b)$.
  The Galois field of $p$ elements, where $p$ is a prime, is denoted by
  $\GF(p)$. If $m$ and $n$ are integers, by $\rest(m,n)$ we shall denote
  the remainder obtained from division of $m$ by $n$.
  \bigskip

  \centerline{\bf 2. The left factorial function}
  \medskip

  \DJ. Kurepa defined in [{\bf Ku71}] an arithmetic function $K(n)$ that he denoted
  by $!n$ and called it the left factorial, by
  $$
       K(n)=\, !n = \sum_{i=0}^{n-1} i!, \quad n\in N^+.
  $$
  In the same paper, Kurepa asked if
  $$
       (!n,n!)= 2, \quad n=2,3,\ldots \qquad .               \tag KH
  $$
  This conjecture, known as the left factorial hypothesis,
  is still an open problem in number theory. There are
  several results which support the truth of the hypothesis. Kurepa
  showed in [{\bf Ku71}] that there are infinitely many $n\in N$ for which KH is
  true. Also, the conjecture is verified by use of computers (Slavi\' c
  for $n<1500$, Wagstaff for $n<50000$, Mijajlovi\'c for $n<310000$, and
  Gogi\'c for $n<1000000$). It is interesting that Kurepa announced the
  positive solution of the problem in 1992, but he never published a
  proof. R.~Guy informed in a letter \v{Z}. Mijajlovi\' c that R.~Bond
  announced a proof of the conjecture too, but the proof was never published.
  The first mention of the left factorial function appeared in [{\bf Ku64}],
  where this function is defined for infinite cardinal numbers as well.
  \medskip

    {\bf 2.1 Some equivalents to KH.}\quad
    There are several statements equivalent to KH. Probably the most
  natural one is the following assertion, which also belongs to Kurepa
  [{\bf Ku71}]:
  $$
     \forall n > 2 \enskip !n\not\equiv 0\enskip (\mmod n).
  $$
  This formulation of the left hypothesis appears in [{\bf Gu}] as problem B44,
  and we shall call this statement also KH. It is not difficult to see that
  this form of KH can be reduced to primes (see [{\bf Ku71}]), i.e. KH is
  equivalent to
  $$
      \forall p\in P\enskip  p>2 \Rightarrow (!p,p)=1         \tag PH
  $$
  where $P$ denotes the set of all primes. Namely if PH fails, then KH
  fails with $p= n$. Conversely, if KH fails, then $n|!n$ for some
  $n>2$. Then there is a prime $p>2$ such that $p|n$ and $p\leq n$.
  If $p=n$, then PH trivially fails. If $p<n$, then
  $$
       !n=\,  !p +\, p! +\ldots+ (n-1)!.
  $$
  Now $p|n$ and $n|!n$ imply $p|\, !n$, and therefore it follows from the
  above relation that $p|\, !p$, contradicting PH. This establishes the
  equivalence of KH and PH.


      If $p$ is a prime, then it is not difficult to establish in $\GF(p)$
   the following identities (see [{\bf Mi}]):
   $$ \align
      !p &= \sum^{p-1}_{k=0} (-1)^{k+1}/k!,
                                 \tag 2.1.1 \\
      !p &= \sum^{p-1}_{k=0} (-1)^k (k+1)(k+2)\dots (p-1).
                                 \tag 2.1.2
      \endalign
   $$
   Since the identity ${p-1\choose k}= (-1)^k$ also holds in $\GF(p)$,
   by (2.1.1) and (2.1.2) the following identities  are true in $\GF(p)$:
   $$ \align
      -!p &= \sum^{p-1}_{k=0} {1\over k!}{p-1\choose k},
                                       \\
      !p &= \sum^{p-1}_{k=0} {p-1\choose k} (k+1)(k+2)\dots (p-1).
      \endalign
   $$
   Therefore, we obtain the following
   \smallskip

   {\caps Theorem 2.1}. KH {\it is equivalent to any of the following
   statements:

   \item{\rm 1.}  For all primes $p$,\enskip $\GF(p) \models
          \displaystyle \sum^{p-1}_{k=0} (-1)^k (k+1)(k+2)\ldots
          (p-1) \not= 0. $
   \smallskip

   \item{\rm 2.}  For all primes $p$, \enskip
          $\displaystyle \sum^{p-1}_{k=0} (-1)^k (k+1)(k+2)\ldots
           (p-1) \not\equiv 0\enskip (\mmod p)$.
   \smallskip

   \item{\rm 3.}  For all primes $p$, \enskip $\GF(p) \models
          \displaystyle \sum^{p-1}_{k=0} \binom{p-1}k (k+1)(k+2)\ldots
          (p-1) \not= 0. $
   \smallskip

   \item{\rm 4.}  For all primes $p$, \enskip
          $\displaystyle \sum^{p-1}_{k=0} \binom{p-1}k (k+1)(k+2)\ldots
          (p-1) \not\equiv 0\enskip (\mmod p)$.
   \smallskip

   \item{\rm 5.}  For all primes $p$, \enskip $\GF(p) \models
          \displaystyle \sum^{p-1}_{k=0} {1\over k!}(-1)^k!
          \not= 0.$
   \smallskip

   \item{\rm 6.}  For all primes $p$, \enskip $\GF(p) \models
          \displaystyle \sum^{p-1}_{k=0} {1\over k!}\binom{p-1}k \neq 0.$
   }
   \medskip

   Here $\GF(p)\models\ldots $ means: in $\GF(p)$ we have $\ldots$.

   The second statement in the above theorem is proved in fact also in
   [{\bf StZi}, (Lemma 2.6)]. There are some other equivalences. In
   [{\bf \v Sa}]
   the following equivalence to KH was proved:
   $$
       \forall n > 2\enskip \left(\sum_{k=2}^{n-1} !k,!n\right)= 2,
   $$
   while in [{\bf St}] KH was proved to be equivalent to
   $$
      \sum_{k=2}^n (k-1)\cdot k! \not\equiv  0\enskip (\mmod n), \quad n>2.
   $$

   {\bf 2.2 Some formulas involving {\rm KH}.}\quad
   There are a number of identities involving $!n$ obtained in [{\bf St}],
   [{\bf StZi}] and [{\bf Ca}].
   Stankovi\'c and \v{Z}i\v{z}ovi\'c (cf. [{\bf St}] and [{\bf StZi}]) proved
   the following identities (we assumed that $K(0)=0$):
   $$ \align
    \sum_{i=0}^n K(i)        &= nK(n-1) + 1,         \quad n\geq 1,
                                  \tag 2.2.1\\
    2\sum_{i=0}^{n-1}iK(i)   &= K(n) + n(n-1)K(n-2), \quad n\geq 2,
                                  \tag 2.2.2\\
    6\sum_{i=0}^{n-1}i^2K(i) &= (2n-1)K(n) + (2n^2-n-2)K(n-2)
                    + 2\cdot n! - 4,     \quad n\geq 2.
                                  \tag 2.2.3\\
      \endalign
   $$

   In connection with these identities, Carlitz (cf. [{\bf Ca}]) considered the
   following sums:
   $$  \align
     Q_m(n) &= \sum_{k=0}^{n-1} k^m K(k), \quad m=0, 1, 2,\ldots,       \\
     R_m(n) &= \sum_{k=0}^{n-1} \binom km K(k).                         \\
       \endalign
   $$
   In the same paper, he proved the following generalizations of (2.2.1--3):
   $$  \align
    R_m(n) &= \binom n{m+1} K(n) -  K_m(n) - K_{m+1}(n),     \tag 2.2.4\\
        &\kern -1cm \text{where}\quad K_m(n)=
                     \sum_{k=0}^{n-1} \binom km k!,    \\
    R_m(n) &= \binom n{m+1} K(n) -
          \sum_{j=0}^m (-1)^{m-j}\binom mj {K(n+j+1)-K(j+1)\over j+1},
                                 \tag 2.2.5\\
    Q_m(n) &= \sum_{k=0}^m k! S(m,k)R_k(n),                  \tag 2.2.6\\
      \endalign
   $$
   where $S(m,k)$ are Stirling numbers of the second kind. Let us note that
   by use of $s(m,k)$, i.e. Stirling numbers of the first kind, we can
   obtain the dual of the identity 2.2.6., that is, we can express
   $R_m(n)$ by $Q_m(n)$. Namely, it is well known that the matrices
   $|| S(m,k) ||$ and $|| s(m,k) ||$ are mutually inverse, therefore,
   from (2.2.6) it follows at once that
   $$
      R_m(n)= {1\over m!} \sum_{k=0}^m s(m,k)Q_k(n).
   $$
   \medskip

   {\bf 2.3 Number theoretical hypotheses related to {\rm KH}.}\quad
   The hypothesis on the alternating factorial stated as the problem B43
   in Guy's monograph [{\bf Gu}] on unsolved problems in number theory
   is similar to KH (stated in [{\bf Gu}] as problem B44).
   Here is the formulation of this problem:

   {\it Let
   $$
      A_n= (n-1)! - (n-2)! + (n-3)! - \ldots + (-1)^n\cdot 1!, \quad n=2,3\ldots
   $$
   Are there infinitely many numbers $n$ such that $A_n$ is a prime?}

   In [{\bf Gu}] it is observed that if there is $n\in N^+$ such that $n+1$ divides
   $A_n$, then $n+1$ will divide $A_m$ for all $m>n$, and there would be
   only finitely many number of prime values of the sequence $A_n$.
   Wagstaff verified   this fact for $n<46340$, while Gogi\'c extended
   this result in his master thesis [{\bf Go}] to $n<1000000$.

      In his paper [{\bf Ku74}], Kurepa asked several question concerning KH.
   He introduced there the statement $H_4(s)$ in the following way:
   $$
       (n\geq 2\wedge s\geq 1) \Rightarrow (K(n),K(n+s))= 2,\quad n,s\in N^+.
                            \tag $H_4(s)$
   $$
   Then Kurepa asked ([{\bf Ku74}], Problem 2.9.) if KH implies $H_4(s)$ for
   all $s\in N^+$. We note that this implication does not hold since,
   for example:
   $$ \matrix \format \r&\enskip\c&\enskip\r&\enskip\l                     \\
      K(7)  &=             &874 &= 2\cdot 19 \cdot 23                      \\
      K(12) &=        &43954714 &= 2\cdot 19 \cdot 31 \cdot 37313          \\
      K(16) &=  &1401602636314 &= 2\cdot 19\cdot 41\cdot 491\cdot 1832213  \\
      K(25) &= &647478071469567844940314 &=
          2\cdot 41\cdot 103 \cdot 2875688099\cdot 26658285041.    \\
    \endmatrix
   $$
   The same examples also show that the strong left factorial hypothesis
   does not hold, as Kurepa formulated it in [{\bf Ku74}]:

   {\it The numbers $K(n)/2$, $n=2, 3, \ldots$ are pairwise relatively prime}.

   In the same paper, Kurepa introduced the  sequence of sets
   $$
     A(r)= \{n\in N^+\, |\, r<n, K(n)\equiv r\enskip (\mmod n) \}.
   $$
   He asked there for a
   description of these sets, and in particular is there any $r$ for
   which $A(r)$ is finite. He also asked if $A(3)= \emptyset $.
   We note here that $467\in A(3)$.
   \medskip

   Kurepa asked in [{\bf Ku71}]  if $!n$ is square-free, with the only exception
   $!3 = 2^2$. This hypothesis, which we shall call KH2, is verified  in [{\bf Mi}]
   for $n\leq 40$ by finding prime decompositions of $!n$ for $n<40$.
   There is a  simple connection between KH and KH2. Namely, if $p$ is a
   prime and $n\geq p$, then $p^2|!n$ implies $p|!n$, and so $p|!p$.
   Hence we obtain
   \smallskip

   {\caps Proposition 2.3.1} {\it {\rm KH} implies that for any $m>1$ there
   are at most finitely many $n$ such that $m^2| !n$.}
   \medskip

   {\bf 2.4 Computational verification of KH.}\quad
      There are simple recurrent formulas for the remainder of $!n$ divided by
   $n$. Using these formulas it easy to check KH and to perform the related
   computation. Let $r_n$ be the sequence defined by $r_n= \rest(!n,n)$,
   $n\in N^+$. The following proposition enables one to design an
   algorithm for computing the values of $r_n$ (cf. [{\bf Mi}, Lemma 2.1-3]):
   \smallskip

   {\caps Proposition 2.4.1} {\it Let $q$ be a prime, and let the finite
   sequences   $s_i$, $t_i$, $v_i$ be defined in $\GF(q)$ in the following
   way:}
$$ \alignat4
 &1. &\enskip s_{q-1} &= 0, &\quad s_i &= 1+is_{i+1},        &\quad i &=q-2,q-3,\ldots, 1.  \\
 &2. &\enskip t_1 &= 0,     &\quad t_i &= (-1)^i + it_{i-1}, &\quad i &=2, 3, \ldots, q-1.  \\
 &3. &\enskip v_1 &= 0,     &\quad v_i &= 1 - i v_{i-1},     &\quad i &= 2,3,\ldots, q-1.   \endalignat
$$

   {\it Then\/  $r_q= s_1 = t_{q-1} = v_{q-1}$}.
   \smallskip

   Observe that $s_q$ is defined by the regressive induction. Using these
   formulas it is easy to develop a simple computer program for verifying
   KH by computing $r_n$. Let KH$(x)$ denote the truth of the left
   factorial hypothesis for all positive integers $n\leq x$.
   Mijajlovi\'c [{\bf Mi}] verified KH$(311009)$ and Gogi\'c [{\bf Go}]
   extended it to all $n<1000000$.

   By simple modification of the above formulas one can obtain in the
   ring ${\bold Z}_{p^2}$, where $p$ is a prime, the following recurrent
   formulas:
   $$ \align
       s_{n-1} &= n,                                      \tag 2.4.1 \\
       s_i &= 1+is_{i+1}, \quad i= n-2, n-3, \ldots, 1,          \\
      \endalign
   $$
   so that $\rest(!n,p^2)=s_1$. Thus using 2.4.1, and assuming KH, in [{\bf Mi}] it
   was proved:  if $m^2|!n$ then $m\geq 1227$.

   By inspection, we see that the total number of arithmetical operations
   used in the verification of KH$(x)$ is
   $$
    A(x)= \sum_{p\leq x} 4p,                               \tag 2.4.2
   $$
   where $p$ in the sum runs over primes.          

   Using the prime number theorem in the form
   $$
       \pi(x)= \sum_{p\leq x} 1 = {x \over \ln x} +
               O\left({x\over \ln^2 x}\right),
   $$
   and integration by parts, we obtain
   $$
      A(x)= {2x^2\over \ln x} + O\left({x^2\over \ln^2 x}\right).
   $$
   Therefore, the growth
   of the number of arithmetical operations used in the verification of
   KH$(x)$ is
   $$
      a(k,x)= {A(kx)\over A(x)} \sim k^2 \quad \text{as}\enskip  x\to \infty.
   $$
   This means, as it was explained in [{\bf Mi}], that the efficiency
   in the verification of KH$(x)$ by use of parallel computers with
   $k$ parallel processors is $\sqrt k$.

   \bigskip

   {\bf 2.5 Left factorial function in complex domain.}\quad
   The gamma-function $\Gamma(z)$ is defined by
   $$
      \Gamma(z)= \int_0^{\infty}e^{-t}t^{z-1}dt  \quad (\Re\, z >0),
                              \tag 2.5.1
   $$
   and for other values of the complex variable $z$ by analytic continuation,
   furnished by the functional equation
   $$
      z\Gamma(z)= \Gamma(z+1).                            \tag 2.5.2
   $$
   Since $\Gamma(n+1)= n!$ for $n\in N$, it follows that
   $$
      K(n)= \sum_{i=0}^{n-1} \Gamma(i+1)=
            \int_0^{\infty}e^{-t}\sum_{i=0}^{n-1}t^idt=
            \int_0^{\infty}e^{-t}{t^n-1\over t-1}dt\quad (n\in N^+).
   $$
   Hence for $\Re z> 0$ it makes sense to define
   $$
      K(z)= \int_0^{\infty}e^{-t}{t^z-1\over t-1}dt,      \tag 2.5.3
   $$
   and since one easily obtains
   $$
      K(z)=  K(z+1)- \Gamma(z+1),                         \tag 2.5.4
   $$
   then (2.5.4) provides analytic continuation of $K(z)$ to the whole
   complex plane. In particular, since $K(1)= \Gamma(1)=1$, it follows
   that $K(0)=0$. Kurepa [{\bf Ku71}] defined $K(z)$ for arbitrary complex
   $z$ by (2.5.3) and (2.5.4). In [{\bf Ku73}] he established that
   $K(z)$ is a meromorphic function having only simple poles at the
   points $z= -1, -3, -4, -5,\ldots$. The residue of $K(z)$ at $z=-1$
   equals to $-1$, and at $z= -n$ ($n=3,4,5,\ldots$) it equals
   $\sum_{k=2}^{n-1}(-1)^{k-1}/k!$. This follows from (2.5.2), (2.5.4) and
   the fact that $\Gamma(z)$ is a meromorphic function with residues
   $(-1)^n/n!$ at simple poles $z=-n$ $(n\in N)$. Kurepa [{\bf Ku73}]
   also studied the zeros of $K(z)$, and showed the asymptotic relations
   $$
     \lim_{x\to\infty} {K(x)\over\Gamma(x)}= 1, \quad
     \lim_{x\to\infty} {K(x)\over\Gamma(x)}= 0,
   $$
   of which the second is a corollary of the first in view of (2.5.2).
   Further results on $K(z)$ as a function of the complex variable $z$
   were obtained by Slavi\' c [{\bf Sl}]. His main result is that
   $$
      K(z)= -{\pi\over e}\text{cotg} \pi z +
        {1\over e}\left(\sum_{k=1}^{\infty}{1\over n!n} + C \right) +
        \sum_{k=0}^{\infty}\Gamma(z-n)
   $$
   holds for all complex $z$, where
   $$
      C= -\int_0^{\infty} e^{-x}\ln x dx = 0.577215\ldots
   $$
   is Euler's constant.

       Formula (2.5.3) is useful for many purposes. For example, for
   $p\geq 3$ it gives
   $$
      !p= K(p)= \int_0^{\infty} e^{-t}{((t-1)+1)^p-1\over t-1}dt\equiv
    \int_0^\infty e^{-t}(t-1)^{p-1} dt\enskip (\mmod p),      \tag 2.5.5
   $$
   since when we expand $((t-1)+1)^p$ by the binomial theorem we can use that
   $\binom pk$ $(k=1,\ldots, p-1)$ is divisible by $p$. But
   $$ \eqalign{
    \int_0^{\infty}e^{-t}(t-1)^{p-1}dt
         &= \int_0^1e^{-t}(t-1)^{p-1}dt +
           {1\over e}\int_0^\infty e^{-u} u^{p-1} du     \cr
         &= \int_0^1 e^{-t}(t-1)^{p-1}dt + (p-1)!/e.         \cr
          }                                                \tag 2.5.6
   $$
   Since the first integral in (2.5.6) is a natural number and
   $$
      0 < \int_0^1 e^{-t}(t-1)^{p-1} dt \leq 1- {1\over e},
   $$
   it follows from (2.5.5) and (2.5.6) that
   $$
     !p \equiv \left[(p-1)!\over e\right]+1 \enskip(\mmod p),  \tag 2.5.7
   $$
   where $[x]$ denotes the integer part of $x$. Therefore, from (2.5.7)
   we can obtain, in view of PH, another equivalent of KH, namely
   $$
     \left[(p-1)!/e\right]\not\equiv -1 \enskip(\mmod p) \quad
                              \text{for}\quad p>2.
   $$
   In connection with (2.5.7) one can define $R(p)$ to be the least
   nonnegative residue of $!p\enskip (\mmod p)$. The evaluation of $R(p)$ is
   rather involved, but perhaps one could try to evaluate the
   summatory function of $R(p)$. The following problem seems to be of
   interest: does there exist a constant $C>0$ such that
   $$
       \sum_{p\leq x} R(p) \sim {Cx^2\over \ln x}\quad   (x\to\infty)?
   $$
   \medskip

   \centerline{\bf 3. Other hypotheses}
   \medskip

    Kurepa presented several problems in number theory at the Problem
    Session of the 5th Balkan Mathematical Congress, held in Belgrade,
    1974. These problems are published as a supplement to [{\bf Ku74}].
    The first problem on this list concern the set
    $$
       P(n)= \{ x\in N^+\colon \{x-2n, x, x+2n\}\subseteq P\}, \quad n\in N^+
                                \tag 3.1
    $$
    where $P$ is the set of prime numbers. Kurepa asked what are the
    properties of $P(n)$, and in particular:
    \smallskip

    \item{P1.} Is P(1)= \{5\}?
    \item{P2.} Is there some $n\in N^+$ such that $P(n)= \emptyset$?
    \smallskip

    We note the following properties of the sequence $P(n)$. The set $P(n)$
    is related to a part of Problem A6 in [{\bf Gu}]. Namely, as noted there,
    it is not known whether
    there are infinitely many sets of three consecutive primes in an
    arithmetic progression, but S. Chowla has shown [{\bf Ch}] this without the
    restriction to consecutive primes. Thus, as $x-2n$, $x$, $x+2n$ is an
    arithmetic progression, we have
    $$
    \bigcup_{n\in N^+} P(n) \quad \text{is infinite}.         \tag 3.2
    $$
    Further, assume $n=3k+1$, $k\in N$ . Then $x-2n\equiv x+1\enskip (\mmod 3)$, and
    $x+2n\equiv x+2\enskip (\mmod 3)$,
    thus $3$ divides $(x-2n)x(x+2n)$. If $x\in P(n)$ then $x-2n=3$, i.e.
    $x= 2n+3$. Hence $P(n)= \emptyset$, or $P(n)$ is an one-element set, i.e.
    $P(n)= \{6k+5\}$, where $6k+5, 12k+7\in P$. For example,
    $P(1)$, $P(4)$, $P(7)$, $P(10)$ are one-element sets, while
    $P(13)= \emptyset$, and this answers questions P1 and P2.

    By (3.2), without any restriction on $n$, there are infinitely many $n$
    such that $P(n)$ is an one-element set.
    However, we may ask if there are infinitely many $k\in N$
    such that $P(n)$ is an one-element set, where $n=3k+1$. We
    already observed that this is the case iff $6k+5, 12k+7 \in P$. We do not
    know the right answer, and obviously this question is related to the
    twin primes conjecture, and to the conjectures 5 and 4 in [{\bf Sh}],
    which in turn would imply that there are infinitely many Mersenne primes.
    On the other hand, as the functions $6k+5$ and $12k+7$ are linearly
    independent, from the Bateman-Horn conjecture [{\bf BaHo}] it would follow
    that the number of $k\leq m$  such that $6k+5, 12k+7\in P$ is asymptotic
    to $\displaystyle C \int_2^m dx/(\log x)^2$, where $C$ is a positive
    constant. In particular, it would follow that there are infinitely many
    $n\in 3N+1$ such that $P(n)$ is an one-element set.

    If $n= 3k+2$, where $k\in N$, we have a similar  conclusion, i.e. $3$ divides
    $(x-2n)x(x+2n)$, and so if $x\in P(n)$ then $x= 2n+3$, thus
    $P(n)= \{6k+7\}$, where $6k+7, 12k+11\in P$. We have also
    a similar discussion as in the case $n=3k+1$.

    Finally, if $n=3k$, $k\in N$, then the problem whether $P(n)$ is infinite
    reduces to the question whether there are infinitely many prime triplets
    $x-6k$, $x$, $x+6k$. This question is related to Problem A9 in [{\bf Gu}]
    and according to the discussion supplemented to the problem,
    it is likely that there are infinitely many such triplets.

    The second problem Kurepa stated in his list concerns the sequence
    $s_n= p_n^2 - p_{n-1} - p_{n+1}$, where $p_n$ is the n-th prime.
    Kurepa asked what could be the sign of the elements of this sequence.
    We note the following:
    \smallskip

    {\caps Lemma.} \quad $p_n^2 > p_{n-1} + p_n + p_{n+1}$ if $p_n\geq 5$.
    \smallskip

    {\it Proof} \, If $p_n\geq 5$, then
    $$
       p_n^2 - p_{n-1} - p_{n+1} \geq 5 p_n -2 p_{n+1} =
           p_n + 2(2p_n- p_{n+1}).
    $$
   By Bertrand's postulate, which says that for every positive integers $m$
   there is a prime in the interval $[m,2m]$, we have $p_{n+1}< 2p_n$, so
   $$ \align
     &p_n^2 - p_{n-1} - p_{n+1} > p_n, \quad i.e.                      \\
     &p_n^2 > p_{n-1} + p_n + p_{n+1}. \qquad\qquad\qquad \diamondsuit \\
      \endalign
   $$
   By the above lemma, we see that for all $p_n\geq 3$,
   $p_n^2 - p_{n-1} - p_{n+1} > 0$. Actually we can show that
   $$
       p_n^2 > \sum_{k=1}^{n+1} p_k\quad  (n\geq n_0).   \tag 3.3
   $$
   Namely, from the prime number theorem it follows that
   $$
        p_n= n(\ln n + O(\ln\ln n)).                \tag 3.4
   $$
   Using (3.4) it follows that the left-hand side of (3.3) is asymptotic
   to $n^2\ln^2 n$, while the right-hand side is
   $$
     \sum_{k=1}^{n+1} k\ln k + O(n^2(\ln\ln n)^2)=
              {n^2\over 2}\ln n + O(n^2(\ln\ln n)^2).
   $$

   In the third problem of his list Kurepa considered the sequence defined by
   $\pi_n= p_n^2 - p_{n-1}p_{n+1}$, and asked what could be the sign
   of members of this sequence, and how often they take the same sign.
   First, let us note that obviously $\pi_n<0$ or $\pi_n>0$. Further,
   this question is related to Problem A14 in [{\bf Gu}]. Namely,
   Erd\H os and Straus call the prime $p_n$ {\it good\/} if
   $p_n^2> p_{n-i}p_{n+i}$ for all $1 \leq i \leq n-1$. Pomerance [{\bf Po}]
   proved that
   there are infinitely many good primes, and therefore there are infinitely
   many $n$ such that $\pi_n>0$.
   Pomerance also proved that
   $$
       \limsup_{n\to\infty} (p_n^2-M(n))= +\infty, \quad \text{where} \quad
       M(n)= \max_{0<i<n} p_{n-i}p_{n+i}.
   $$
   Now suppose
   that $p_{n-1}$ and $p_n= p_{n-1}+2$ are twin primes. Then
   $p_{n+1}\geq p_n + 6$, thus (if $p_n \geq 3$)
   $$
       p_n^2 - p_{n-1} p_{n+1} \leq
           p_{n-1}^2 +4p_{n-1} +4 - p_{n-1}(p_{n-1}+6)=
           -2p_{n-1} +4 <0.
   $$
   Hence, from the twin prime conjecture it would follow that there are
   infinitely many $n\in N^+$ such that $\pi_n< 0$. For some further
   inequalities involving $p_n$ we refer the reader to the monograph
   of Mitrinovi\'c and Popadi\'c [{\bf MiPo}].

    The last problem in number theory (Problem 4) from the Kurepa's list
    concerns the left factorial hypothesis, and we discussed it already
    in the previous section.

\vfill\eject


\bsk\ssk

\begingroup
\frenchspacing
\eightpoint\rm
\cent{REFERENCES}\nobreak
\msk\ssk\nobreak
\setbox1=\hbox{[{\bf2}]\enskip}\parindent=\wd1

   \item{[BaHo]} P. Bateman, R.Horn, {\it A heuristic asymptotic formula
      concerning the distibution of prime numbers}, Math. Comp. {\bf 16}
      (1962), 363--367.
   \item{[Ca]} L. Carlitz, {\it A note on the left factorial function},
      Math. Balkan. {\bf 5:6} (1975), 37--42.
   \item{[Ch]} S. Chowla, {\it There exists an infinity of\/
      $3$-combinations of primes in A.P.}, Proc. Lahore, Philos. Soc. {\bf 6}
      (2) (1944), 15--16.
   \item{[Go]} G. Gogi\'c, {\it Parallel Algorithms in Arithmetic},
      Master thesis, Belgrade University, 1991, (in Serbian).
   \item{[Gu]} R. Guy, {\it Unsolved Problems in Number Theory},
      Springer-Verlag, 1981.
   \item{[Ku64]} \DJ. Kurepa, {\it Factorials of cardinal numbers and
      trees}, Glasnik Mat. Fiz. Astr. {\bf 19} (1-2) (1964), 7--21.
   \item{[Ku71]} \DJ. Kurepa, {\it On the left factorial function $n!$},
      Math. Balkan. {\bf 1} (1971), 147--153.
   \item{[Ku73]} \DJ. Kurepa {\it Left factorial function in complex domain},
      Math. Balkan. {\bf 3} (1973), 297--307.
   \item{[Ku74]} \DJ. Kurepa, {\it On some new left factorial propositions},
      Math. Balkan. {\bf 4} (1974), 383--386.
   \item{[Mi]} \v{Z}. Mijajlovi\'c, {\it On some formulas involving $!n$ and the
      verification of the $!n$-hypothesis by use of computers}, Publ. Inst.
      Math., {\bf 47}({\bf 61}) (1990), 24--32.
   \item{[MiPo]} D.S. Mitrinovi\'c, M.S. Popadi\'c, {\it Inequalities
      in Number Theory}, University of Ni\v s, Ni\v s, 1978.
   \item{[Na]} W. Narkiewicz, {\it Classical Problems in Number Theory},
      Monograf mat., {\bf 62}, PWN, Warszawa, 1986.
   \item{[Po]} C. Pomerance, {\it The prime number graph}, Math. Comp.
      {\bf 33} (1979), 399--408.
   \item{[Ri]} H. Riesel, {\it Prime Numbers and Computer Methods for
      Factorization}, Birk\-h\"aus\-er, 1985.
   \item{[Sh]} D. Shanks, {\it Solved and Unsolved Problems in Number
      theory}, 3ed., Chelsea, New York, 1985.
   \item{[Sl]} D.V. Slavi\'c, {\it On the left factorial function of the
      complex argument}, Math. Balkan. {\bf 3} (1973), 472--477.
   \item{[St]} J. Stankovi\'c, {\it \"Uber einige Relationen zwischen
      Fakult\"aten und den linken Fakult\"aten}, Math. Balkan.
      {\bf 3} (1973), 488--497.
   \item{[StZi]} J. Stankovi\'c, M. \v{Z}i\v{z}ovi\'c, {\it Noch einige
      Relationen zwischen den Fakult\"aten und den linken Fakult\"aten},
      Math. Balkan. {\bf 4} (1974), 555--559.
   \item{[\v Sa]} Z. \v Sami, {\it On the M-hypothesis of D. Kurepa}, Ibidem,
      530--532.
   \item{[Vu]} V. Vuleti\'c, {\it Tabulation of the functions}:
      $\Gamma(n+1)=n!$, $K(n)= !n$, $r(n)$, ${\nu}_s(n)=\nu(s,n)$,
      Math. Balkan. {\bf 4} (1974), 675-706.


   \bsk
   \parindent=25pt\parskip0pt

\vskip2cm
     Matemati\v{c}ki Institut     \hfill (Received 29 06 1994)

     Knez Mihailova 35

     11000 Beograd, P.P. 367

     (Serbia) Yugoslavia

     email: {\tt  aivic\@matf.bg.ac.yu} (A. Ivi\'c)

     {\tt emijajlo\@ubbg.etf.bg.ac.yu} (\v Z. Mijajlovi\'c)

\vskip 1cm

\parindent=5cm
   \bye